\catcode`\@=11

\thicklines
\newskip\Einheit \Einheit=.6cm
\newcount\xcoord \newcount\ycoord
\newdimen\xdim \newdimen\ydim \newdimen\PfadD@cke \newdimen\Pfadd@cke
\def\PfadDicke#1{\PfadD@cke#1 \divide\PfadD@cke by2 
\Pfadd@cke\PfadD@cke \multiply\PfadD@cke by2}
\long\def\LOOP#1\REPEAT{\def\BODY{#1}\ITERATE}
\def\ITERATE{\BODY \let\next\ITERATE \else\let\next\relax\fi \next}
\let\REPEAT=\fi
\def\Punkt{\hbox{\raise-2pt\hbox to0pt{\hss\scriptsize$\bullet$\hss}}}

\def\DuennPunkt(#1,#2){\unskip
  \raise#2 \Einheit\hbox to0pt{\hskip#1 \Einheit
          \raise-1.5pt\hbox to0pt{\hss\tiny$\bullet$\hss}\hss}}
		  
\def\NormalPunkt(#1,#2){\unskip
  \raise#2 \Einheit\hbox to0pt{\hskip#1 \Einheit
          \raise-3pt\hbox to0pt{\hss\large$\bullet$\hss}\hss}}
\def\DickPunkt(#1,#2){\unskip
  \raise#2 \Einheit\hbox to0pt{\hskip#1 \Einheit
          \raise-4pt\hbox to0pt{\hss\Large$\bullet$\hss}\hss}}
\def\Kreis(#1,#2){\unskip
  \raise#2 \Einheit\hbox to0pt{\hskip#1 \Einheit
          \raise-4pt\hbox to0pt{\hss\Large$\circ$\hss}\hss}}
\def\Diagonale(#1,#2)#3{\unskip\leavevmode
  \xcoord#1\relax \ycoord#2\relax
      \raise\ycoord \Einheit\hbox to0pt{\hskip\xcoord \Einheit
         \unitlength\Einheit
         \line(1,1){#3}\hss}}
\def\AntiDiagonale(#1,#2)#3{\unskip\leavevmode
  \xcoord#1\relax \ycoord#2\relax \advance\xcoord by -0.05\relax
      \raise\ycoord \Einheit\hbox to0pt{\hskip\xcoord \Einheit
         \unitlength\Einheit
         \line(1,-1){#3}\hss}}
\def\Pfad(#1,#2),#3\endPfad{\unskip\leavevmode
  \xcoord#1 \ycoord#2 \thicklines\ZeichnePfad#3\endPfad\thinlines}
\def\ZeichnePfad#1{\ifx#1\endPfad\let\next\relax
  \else\let\next\ZeichnePfad
    \ifnum#1=1
      \raise\ycoord \Einheit\hbox to0pt{\hskip\xcoord \Einheit
         \vrule height\Pfadd@cke width1 \Einheit depth\Pfadd@cke\hss}%
      \advance\xcoord by 1
     \else\ifnum#1=2
      \raise\ycoord \Einheit\hbox to0pt{\hskip\xcoord \Einheit
         \unitlength\Einheit
         \line(0,1){1}\hss}
      \advance\xcoord by 0
      \advance\ycoord by 1
 \else\ifnum#1=3
      \raise\ycoord \Einheit\hbox to0pt{\hskip\xcoord \Einheit
         \unitlength\Einheit
         \line(1,1){1}\hss}
      \advance\xcoord by 1
      \advance\ycoord by 1
    \else\ifnum#1=4
      \raise\ycoord \Einheit\hbox to0pt{\hskip\xcoord \Einheit
         \unitlength\Einheit
         \line(1,-1){1}\hss}
      \advance\xcoord by 1
      \advance\ycoord by -1
   \else\ifnum#1=5
      \raise\ycoord \Einheit\hbox to0pt{\hskip\xcoord \Einheit
         \unitlength\Einheit
         \line(2,1){2}\hss}
      \advance\xcoord by 2
      \advance\ycoord by 1
	  \else\ifnum#1=6
      \raise\ycoord \Einheit\hbox to0pt{\hskip\xcoord \Einheit
         \unitlength\Einheit
         \line(2,-1){2}\hss}
      \advance\xcoord by 2
      \advance\ycoord by -1
	  \else\ifnum#1=7
      \raise\ycoord \Einheit\hbox to0pt{\hskip\xcoord \Einheit
         \unitlength\Einheit
         \line(3,1){3}\hss}
      \advance\xcoord by 3
      \advance\ycoord by 1
	  \else\ifnum#1=8
      \raise\ycoord \Einheit\hbox to0pt{\hskip\xcoord \Einheit
         \unitlength\Einheit
         \line(3,-1){3}\hss}
      \advance\xcoord by 3
      \advance\ycoord by -1
    \fi\fi\fi\fi\fi\fi\fi\fi
  \fi\next}
\def\hSSchritt{\leavevmode\raise-.4pt\hbox 
to0pt{\hss.\hss}\hskip.2\Einheit
  \raise-.4pt\hbox to0pt{\hss.\hss}\hskip.2\Einheit
  \raise-.4pt\hbox to0pt{\hss.\hss}\hskip.2\Einheit
  \raise-.4pt\hbox to0pt{\hss.\hss}\hskip.2\Einheit
  \raise-.4pt\hbox to0pt{\hss.\hss}\hskip.2\Einheit}
\def\vSSchritt{\vbox{\baselineskip.2\Einheit\lineskiplimit0pt
\hbox{.}\hbox{.}\hbox{.}\hbox{.}\hbox{.}}}
\def\DSSchritt{\leavevmode\raise-.4pt\hbox to0pt{%
  \hbox to0pt{\hss.\hss}\hskip.2\Einheit
  \raise.2\Einheit\hbox to0pt{\hss.\hss}\hskip.2\Einheit
  \raise.4\Einheit\hbox to0pt{\hss.\hss}\hskip.2\Einheit
  \raise.6\Einheit\hbox to0pt{\hss.\hss}\hskip.2\Einheit
  \raise.8\Einheit\hbox to0pt{\hss.\hss}\hss}}
\def\dSSchritt{\leavevmode\raise-.4pt\hbox to0pt{%
  \hbox to0pt{\hss.\hss}\hskip.2\Einheit
  \raise-.2\Einheit\hbox to0pt{\hss.\hss}\hskip.2\Einheit
  \raise-.4\Einheit\hbox to0pt{\hss.\hss}\hskip.2\Einheit
  \raise-.6\Einheit\hbox to0pt{\hss.\hss}\hskip.2\Einheit
  \raise-.8\Einheit\hbox to0pt{\hss.\hss}\hss}}
\def\SPfad(#1,#2),#3\endSPfad{\unskip\leavevmode
  \xcoord#1 \ycoord#2 \ZeichneSPfad#3\endSPfad}
\def\ZeichneSPfad#1{\ifx#1\endSPfad\let\next\relax
  \else\let\next\ZeichneSPfad
    \ifnum#1=1
      \raise\ycoord \Einheit\hbox to0pt{\hskip\xcoord \Einheit
         \hSSchritt\hss}%
      \advance\xcoord by 1
    \else\ifnum#1=2
      \raise\ycoord \Einheit\hbox to0pt{\hskip\xcoord \Einheit
        \hbox{\hskip-2pt \vSSchritt}\hss}%
      \advance\ycoord by 1
    \else\ifnum#1=3
      \raise\ycoord \Einheit\hbox to0pt{\hskip\xcoord \Einheit
         \DSSchritt\hss}
      \advance\xcoord by 1
      \advance\ycoord by 1
    \else\ifnum#1=4
      \raise\ycoord \Einheit\hbox to0pt{\hskip\xcoord \Einheit
         \dSSchritt\hss}
      \advance\xcoord by 1
      \advance\ycoord by -1
    \fi\fi\fi\fi
  \fi\next}
\def\Koordinatenachsen(#1,#2){\unskip
 \hbox to0pt{\hskip-.5pt\vrule height#2 \Einheit width.5pt depth1 
\Einheit}%
 \hbox to0pt{\hskip-1 \Einheit \xcoord#1 \advance\xcoord by1
    \vrule height0.25pt width\xcoord \Einheit depth0.25pt\hss}}
\def\Koordinatenachsen(#1,#2)(#3,#4){\unskip
 \hbox to0pt{\hskip-.5pt \ycoord-#4 \advance\ycoord by1
    \vrule height#2 \Einheit width.5pt depth\ycoord \Einheit}%
 \hbox to0pt{\hskip-1 \Einheit \hskip#3\Einheit 
    \xcoord#1 \advance\xcoord by1 \advance\xcoord by-#3 
    \vrule height0.25pt width\xcoord \Einheit depth0.25pt\hss}}
\def\Gitter(#1,#2){\unskip \xcoord0 \ycoord0 \leavevmode
  \LOOP\ifnum\ycoord<#2
    \loop\ifnum\xcoord<#1
      \raise\ycoord \Einheit\hbox to0pt{\hskip\xcoord 
\Einheit\Punkt\hss}%
      \advance\xcoord by1
    \repeat
    \xcoord0
    \advance\ycoord by1
  \REPEAT}
\def\Gitter(#1,#2)(#3,#4){\unskip \xcoord#3 \ycoord#4 \leavevmode
  \LOOP\ifnum\ycoord<#2
    \loop\ifnum\xcoord<#1
      \raise\ycoord \Einheit\hbox to0pt{\hskip\xcoord 
\Einheit\Punkt\hss}%
      \advance\xcoord by1
    \repeat
    \xcoord#3
    \advance\ycoord by1
  \REPEAT}
\def\Label#1#2(#3,#4){\unskip \xdim#3 \Einheit \ydim#4 \Einheit
  \def\lo{\advance\xdim by-.5 \Einheit \advance\ydim by.5 \Einheit}%
  \def\llo{\advance\xdim by-.25cm \advance\ydim by.5 \Einheit}%
  \def\loo{\advance\xdim by-.5 \Einheit \advance\ydim by.25cm}%
  \def\o{\advance\ydim by.25cm}%
  \def\ro{\advance\xdim by.5 \Einheit \advance\ydim by.5 \Einheit}%
  \def\rro{\advance\xdim by.25cm \advance\ydim by.5 \Einheit}%
  \def\roo{\advance\xdim by.5 \Einheit \advance\ydim by.25cm}%
  \def\l{\advance\xdim by-.30cm}%
  \def\r{\advance\xdim by.30cm}%
  \def\lu{\advance\xdim by-.5 \Einheit \advance\ydim by-.6 \Einheit}%
  \def\llu{\advance\xdim by-.25cm \advance\ydim by-.6 \Einheit}%
  \def\luu{\advance\xdim by-.5 \Einheit \advance\ydim by-.30cm}%
  \def\u{\advance\ydim by-.30cm}%
  \def\ru{\advance\xdim by.5 \Einheit \advance\ydim by-.6 \Einheit}%
  \def\rru{\advance\xdim by.25cm \advance\ydim by-.6 \Einheit}%
  \def\ruu{\advance\xdim by.5 \Einheit \advance\ydim by-.30cm}%
  #1\raise\ydim\hbox to0pt{\hskip\xdim
     \vbox to0pt{\vss\hbox to0pt{\hss$#2$\hss}\vss}\hss}%
}
\catcode`\@=12

\documentclass[12pt]{article}

\usepackage{amsmath,amsthm}
\usepackage{amssymb}
\usepackage{color}
\usepackage{xspace}
\usepackage[colorlinks=true,
linkcolor=green,
filecolor=brown,
citecolor=green]{hyperref}

\def\blue{\textcolor{blue} }

\begin{document}
\newtheorem{lemma}{Lemma}
\newtheorem{theorem}{Theorem}
\newtheorem{prop}{Proposition}
\newtheorem{cor}{Corollary}
\mbox{ }
\vspace*{-25mm}
\begin{center}
{\Large
Kreweras's Narayana Number Identity Has a Simple Dyck Path 
Interpretation \\ 
}
\vspace{10mm}
DAVID CALLAN  \\
Department of Statistics  \\
\vspace*{-2mm}
University of Wisconsin-Madison  \\
\vspace*{-2mm}
Medical Science Center \\
\vspace*{-2mm}
1300 University Ave  \\
\vspace*{-2mm}
Madison, WI \ 53706-1532  \\
{\bf callan@stat.wisc.edu}  \\

\vspace{2mm}
%Oct 28, 2005
\end{center}

\begin{abstract}
We show that an identity of Kreweras for the Narayana numbers counts Dyck paths with a given number of peaks by number of peak plateaus, where a peak plateau is a run of consecutive peaks that is immediately preceded by an upstep and followed by a downstep.
\end{abstract}

\vspace{8mm} 
Germain Kreweras \cite{kreweras1967,NarayanaLP} proved 
that
\[
\frac{1}{n+r+1}\binom{n+r+1}{r}\binom{n+r+1}{n}=
\sum_{s=0}^{r}\frac{1}{n}\binom{n}{r-s}\binom{n}{r-s+1}\binom{2n+s}{2n}.
\] 
Replacing $r$ by $r-1$, this identity can be rewritten as
\[
\hspace*{40mm}N(n+r,r)=\sum_{s=1}^{r}N(n,s)\binom{2n+r-s}{r-s}\hspace*{40mm} (*)
\]
where $N(n,k):=\frac{1}{n}\binom{n}{k}\binom{n}{k-1}$ is the Narayana 
number, well known to count Dyck $n$-paths containing $k$ peaks 
\cite[Ex.\,6.36, p.\,237]{ec2}.
A Dyck $n$-path is a path of $n$ upsteps $U$ and $n$ downsteps $D$ that never 
dips below the line joining its start and end points.  
A peak is an occurrence of $UD$.
We will show that $(*)$ counts Dyck $(n+r)$-paths with $r$ peaks by 
number $s$ of peak plateaus, where
a peak plateau is a run of consecutive peaks that is preceded by an 
upstep and followed by a downstep, that is, a subpath of the form 
$U(UD)^{i}D$ with $i\ge 1$. For example, the path illustrated below has 5 
peaks but only 2 peak plateaus (in blue, following upsteps 5 and 8 
respectively).

\Einheit=0.8cm
\[
\SPfad(-9,0),111111111111111111\endSPfad
\Pfad(-9,0),3343433\endPfad
\Pfad(2,3),4443\endPfad
\Pfad(8,1),4\endPfad
\blue{\Pfad(-2,3),3434\endPfad
\Pfad(6,1),34\endPfad }
\DuennPunkt(-9,0)
\DuennPunkt(-8,1)
\DuennPunkt(-7,2)
\DuennPunkt(-6,1)
\DuennPunkt(-5,2)
\DuennPunkt(-4,1)
\DuennPunkt(-3,2)
\DuennPunkt(-2,3)
\DuennPunkt(-1,4)
\DuennPunkt(-0,3)
\DuennPunkt(1,4)
\DuennPunkt(2,3)
\DuennPunkt(3,2)
\DuennPunkt(4,1)
\DuennPunkt(5,0)
\DuennPunkt(6,1)
\DuennPunkt(7,2)
\DuennPunkt(8,1)
\DuennPunkt(9,0)
\Label\o{\textrm{{\scriptsize 1}}}(-8.7,0.3)
\Label\o{\textrm{{\scriptsize 2}}}(-7.7,1.3)
\Label\o{\textrm{{\scriptsize 3}}}(-5.7,1.3)
\Label\o{\textrm{{\scriptsize 4}}}(-3.7,1.3)
\Label\o{\textrm{{\scriptsize 5}}}(-2.7,2.3)
\Label\o{\textrm{{\scriptsize 6}}}(-1.7,3.3)
\Label\o{\textrm{{\scriptsize 7}}}(0.3,3.3)
\Label\o{\textrm{{\scriptsize 8}}}(5.3,0.3)
\Label\o{\textrm{{\scriptsize 9}}}(6.3,1.3)
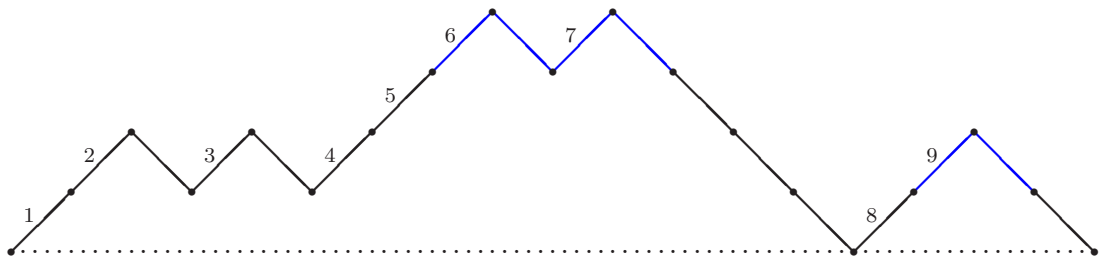
\Label\u{\textrm{{\footnotesize Dyck path with 2 peak 
plateaus}}}(0,-0.3)
\]
%\centerline{{\footnotesize Dyck path with 2 peak plateaus}}
%\vspace*{2mm}
%

\noindent Given a Dyck $(n+r)$-path with $r$ peaks, delete all its peaks. The 
result is a Dyck $n$-path with $s$ peaks where $s$ is the number of peak 
plateaus in the original path. Each Dyck $n$-path with $s$ peaks is 
hit $\binom{2n+r-s}{r-s}$ times: the preimages are found by inserting 
$r-s$ $UD$s, repetition allowed, among the $2n+1$ vertices of the 
$n$-path---$\binom{2n+r-s}{r-s}$ choices---and then inserting 
a $UD$ at each of the $s$ original peak vertices, and $(*)$ follows.

\end{document}